# EQUILIBRIUM MEASURE FOR ONE-DIMENSIONAL LORENZ-LIKE EXPANDING MAPS

BRONZI, M. A. AND OLER, J. G.

**Abstract.** Let $L : [0,1] \setminus \{d\} \to [0,1]$ be a one-dimensional Lorenz like expanding map ($d$ is the point of discontinuity), $\mathcal{P} = \{(0,d), (d,1)\}$ be a partition of $[0,1]$ and $C^\alpha([0,1], \mathcal{P})$ the set of piecewise Hölder-continuous potential of $[0,1]$ with the usual $C^0$ topology. In this context, we prove, improving a result of [2], that piecewise Hölder-continuous potential $\phi$ satisfying $\max\left\{\limsup_{n\to\infty} \frac{1}{n}(S_n\phi)(0), \limsup_{n\to\infty} \frac{1}{n}(S_n\phi)(1)\right\} < P_{\text{top}}(\phi, T)$ support an unique equilibrium state. Indeed, we prove there exists an open and dense subset $\mathcal{H}$ of $C^\alpha([0,1], \mathcal{P})$ such that, if $\phi \in \mathcal{H}$, then $\phi$ admits one equilibrium measure.

## Contents

1. Introduction

A piecewise expanding map $(X, \mathcal{P}, T)$ is a locally connected compact metric space $X$ together with a partition $\mathcal{P}$, which is a finite collection of nonempty pairwise disjoint open subsets of $X$ with dense union, and a map $T : \cup_{P \in \mathcal{P}} P \to X$ such that for each $P \in \mathcal{P}$, each restriction $T|_P$ can be extended to an expanding homeomorphism from a neighborhood of $\overline{P}$ onto one of $\overline{T(P)}$. Here $\overline{P}$ denotes the closure of $P$. An invariant probability measure $\mu$ is called an equilibrium measure for $\phi : M \to \mathbb{R}$ if $h_\mu(L) + \int \phi \, d\mu$ is well-defined and maximal. In this context, Buzzi and Sarig proved in [2] that any piecewise expanding map $T : X \to X$, strongly topologically transitive in the sense that for all nonempty open sets $U$, $T(X) \subset \cup_{k \geq 0} T^k(U)$, with a piecewise Hölder-continuous potential $\phi$ satisfying $P_{\text{top}}(\phi, \partial \mathcal{P}, T) < P_{\text{top}}(\phi, T)$ supports a unique equilibrium state. Here $P_{\text{top}}(\phi, \partial S, T)$ is the topological pressure of a subset $S \subseteq X$ (not necessarily invariant).

In [11] Pesin and Zhang proved that fixed a piecewise expanding map $(X, \mathcal{P}, T)$ the class of potentials, admitting the unique equilibrium measure, is larger than the class of Hölder continuous potentials, i.e., they construct a family of continuous (but not Hölder continuous) potentials $\varphi_c$ exhibiting phase transitions, i.e., there exists a critical value $c_0 > 0$ such that for every $0 < c < c_0$ there is a unique equilibrium measure for $\varphi_c$ which is supported on $(0, 1]$ and for $c < c_0$ the equilibrium measure is the Dirac measure at $0$.

The goal of this paper is to show that if $L : [0, 1] \setminus \{d\} \to [0, 1]$ is a one-dimensional Lorenz like expanding, then the class of potentials piecewise Hölder continuous, admitting the unique equilibrium measure, too is substantially large. More precisely, adapting results of Buzzi and Sarig [2], we construct a open and dense subset $\mathcal{H}$ of $C^\alpha([0, 1], \mathcal{P})$ in $C^0$ topology such that each $\phi \in \mathcal{H}$ admits at most one equilibrium measure. This set $\mathcal{H}$ can be characterized in terms of a regularity condition on the pressure function of the boundary of continuity domains of the dynamics, i.e., if $\phi \in \mathcal{H}$, then $P_{\text{top}}(\phi, T) > \max \left\{ \limsup_{n \to \infty} \frac{1}{n}(S_n\phi)(0), \limsup_{n \to \infty} \frac{1}{n}(S_n\phi)(1) \right\}$ (see Section 3.1).

Acknowledgment

J. O. and M. B. want to thank ICMC-USP for the kind hospitality specially to Ali Tahzibi for suggesting the problem and the many helpful suggestions during the preparation of the paper and wish to express their thanks to Daniel Smania for several helpful comments. They also thank warm hospitality of University Federal de Alagoas (specially to Krerley Oliveira for invitation and conversation).

## 2. Setting and Statements

Lorenz maps originally arise from the study of geometric models for the Lorenz equations ([4], [5], [7], [12], [14]). The Poincaré map has an invariant stable foliation of vertical lines. Since $\mathcal{P}$ takes vertical lines to vertical lines, it induces a map $f : [0,1] \setminus \{d\} \to \mathbb{R}$. This map is an example of an *one-dimensional Lorenz-like expanding map*. The precise definition is as following:

**Definition 2.1.** *A one-dimensional Lorenz-like expanding map is a function $L : [0,1] \to [0,1]$ which has the following properties:*

(L.1) *$L$ has a unique discontinuity at $x = d$ and*
$$L(d^+) = \lim_{x \to d^+} L(x) = 0, \ L(d^-) = \lim_{x \to d^-} L(x) = 1;$$

(L.2) *For any $x \in [0,1] \setminus \{d\}$, $L'(x) > \sqrt{2}$ and the lateral limits at $x = d$ are $L'(d^+) = -\infty$, $L'(d^-) = +\infty$;*

(L.3) *Each inverse branch of $L$ extends to a $C^{1+\theta}, \theta > 0$, function over $[L(0), 1]$ or $[0, L(1)]$ and if $g$ denotes any of these inverse branches, $g'(x) \leq \lambda < 1$.*

**Definition 2.2.** *We say that a continuous map $f : M \to M$ is locally eventually onto (LEO) if given an open set $U \subset M$ there exists $k \in \mathbb{N}$ such that $f^k(U) = M$.*

**Remark 2.1.** *In [14] Williams showed that if $L$ is the one-dimensional Lorenz-like expanding map, then $L$ is LEO.*

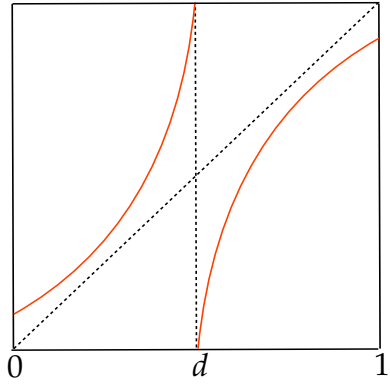

Figure 1. One-dimensional Lorenz-like expanding map.

We denote by $\mathcal{P}$ the natural partition of $[0,1]$, i.e., $\mathcal{P} = \{(0,d),(d,1)\}$. The *boundary* of $\mathcal{P}$ is $\partial \mathcal{P} := \{0, d, 1\}$. Also, define $\mathcal{P}^{(n-1)} = \{P_0 \cap L^{-1}(P_1) \cap \cdots \cap L^{-n+1}(P_{n-1}) \neq \emptyset \mid P_i \in \mathcal{P}\}$. Let $\phi : [0,1] \to \mathbb{R}$ such that $\sup(\phi) < \infty$. This map $\phi$ is a *piecewise Hölder continuous potential* if the restriction to any elements of $\mathcal{P}$ is Hölder continuous, i.e., for all $x, y$ in the same element of $\mathcal{P}$, $|\phi(x) - \phi(y)| \leq K|x-y|^\alpha$ for some $\alpha > 0$, $K < \infty$. Let

$C^\alpha([0,1], \mathcal{P}) := \{\phi : [0,1] \to \mathbb{R} : \phi \text{ is piecewise Hölder continuous potential}\}.$

According Buzzi-Sarig [2], the *pressure of subset* $S \subset [0,1]$ and $\phi \in C^\alpha([0,1], \mathcal{P})$ is defined as

$$P_{top}(\phi, S, L) = \limsup_{n \to \infty} \frac{1}{n} \log \left( \sum_{C_n \in \mathcal{P}^{(n-1)} : S \cap \overline{C_n} \neq \emptyset} \sup_{x \in C_n} e^{S_n \phi(x)} \right),$$

where for $x \in C_n$, $C_n \in \mathcal{P}^{(n-1)}$, the Birkohoff average $S_n \phi(x) = \sum_{j=1}^{n-1} \phi(L^{n-1}(x))$ is well defined. The *topological pressure* of $L$ for $\phi \in C^\alpha([0,1], \mathcal{P})$ is defined by $P(\phi, L) = P_{top}(\phi, [0,1], L)$. In this context, Buzzi and Sarig proved the following theorem:

**Theorem 2.1** ([2])**.** *Let $(X, \mathcal{P}, L)$ be a piecewise expanding map with $\phi \in C^\alpha([0,1], \mathcal{P})$. Assume that $P_{top}(\phi, \partial \mathcal{P}, L) < P_{top}(\phi, L)$. If $\mathcal{M}_L(X)$ denotes the set of invariant measures, then:*

(i) $P_{top}(\phi, L) = \sup_{\mu \in \mathcal{M}_L(X)} \left\{ h_\mu(L) + \int \phi \, d\mu \right\}$ *and this supremum is realized by at least one measure;*
(ii) *there exist at most finitely many ergodic equilibrium measures; and*
(iii) *if, additionally, $L$ is strongly topologically transitive in the sense that for all non-empty open sets $U$, $\cup_{k \geq 0} L^k(U) \supseteq L(X)$, then there exists a unique equilibrium measure.*

In this context, we have can state the main result of the paper:

**Theorem A.** *Let $L : [0,1] \setminus \{d\} \to [0,1]$ be an one-dimensional Lorenz-like expanding map and $C^\alpha([0,1], \mathcal{P})$ the set of piecewise Hölder-continuous potential of [0,1]. Then there exists an open and dense subset $\mathcal{H}$ of $C^\alpha([0,1], \mathcal{P})$ in the $C^0$ topology such that, for every $\phi \in \mathcal{H}$ admit a unique equilibrium measure.*

The idea of the proof of the Theorem A is as following: we use a result of Buzzi-Sarig which gives a sufficient condition for the existence and uniqueness of equilibrium measure under a regularity condition on the pressure of the boundary of continuity domains of the dynamics. We express the pressure of the boundary of the corresponding partition for the one-dimensional Lorenz-like expanding map in terms of the Birkhoff average of the discontinuity and by a small perturbation of the potential along periodic points of sufficiently large period we guarantee the regularity condition of Buzzi-Sarig. To obtain adequate periodic points we use the conjugacy with generalized $\beta$-transformations and the notion of cutting times in one dimensional dynamics (see [1, 6]).

3. Proof of the Theorem A

3.1. **Construction of set $\mathcal{H}$.** As $L$ is not defined in $d \in [0,1]$ we make the following convention: $S_n \phi(d^+)$ ($S_n \phi(d^-)$) is the limit to the right a (resp. to the left) of $z$ of the function $S_n \phi(z)$. More precisely, let $\phi \in C^\alpha([0,1], \mathcal{P})$ be a continuous function, then for any $n \in \mathbb{N}$ we defined

$$S_n\phi(d^\pm) = \lim_{z \to d^\pm} \sum_{i=0}^{n-1} \phi(L^i(z)).$$

By definition $L(d^+) = 0$ and $L(d^-) = 1$, so we conclude that:

$$\limsup_{n\to\infty} \frac{1}{n}S_n\phi(d^+) = \limsup_{n\to\infty} \frac{1}{n}S_n\phi(0) \quad \text{and} \quad \limsup_{n\to\infty} \frac{1}{n}S_n\phi(d^-) = \limsup_{n\to\infty} \frac{1}{n}S_n\phi(1).$$

See Lemma 3.1 to verify that the above relations are well defined.

**Lemma 3.1.** *Let $L : [0,1] \setminus \{d\} \to [0,1]$ be a one-dimensional Lorenz-like expanding map and consider $\phi \in C^\alpha([0,1], \mathcal{P})$.*

  *(i) If does not exist $n_0 \in \mathbb{N}$ such that $L^{n_0}(0) = d$, then*

$$\limsup_{n\to\infty} \frac{1}{n}S_n\phi(d^+) = \limsup_{n\to\infty} \frac{1}{n}S_n\phi(0).$$

  *(ii) If there exists $n_0 \in \mathbb{N}$ such that $L^{n_0}(0) = d$, then*

$$\limsup_{n\to\infty} \frac{1}{n}S_n\phi(d^+) = \frac{1}{n_0}S_{n_0}\phi(0).$$

*The same conclusion holds for $d^-$ replacing 0 for 1.*

*Proof.* The first item comes easily from the definition of Birkhoff sum. Indeed,

$$S_n\phi(d^+) = \lim_{x\to d^+} S_n\phi(x) = \phi(d^+) + \sum_{j=0}^{n-2} \phi(L^j(0)) = \phi(d^+) + S_{n-1}\phi(0),$$

then

$$\limsup_{n\to\infty} \frac{1}{n}S_n\phi(d^+) = \limsup_{n\to\infty} \frac{1}{n}S_n\phi(0).$$

To prove (ii) let $n = Kn_0 + l$, with $l < n_0$. Thus

$$S_n\phi(d^+) = \phi(d^+) + S_{Kn_0+l-1}\phi(0) = \phi(d^+) + KS_{n_0}\phi(0) + S_{l-1}\phi(0),$$

where the second equality in the above equation is because $L$ is piecewise increasing and $\phi$ is continuous. Indeed, for $x$ close enough to $d$ and $x > d$, one has $L^{n_0}(x) > L^{n_0}(0) = d$. Thefore,

$$\frac{1}{n}S_n\phi(d^+) = \frac{\phi(d^+)}{n} + \frac{K}{n}S_{n_0}\phi(0) + \frac{l}{n}S_{l-1}\phi(0).$$

Letting $n \to \infty$ we have:

$$\limsup_{n\to\infty} \frac{1}{n}S_n\phi(d^+) = \frac{1}{n_0}S_{n_0}\phi(0).$$

The same argument gives similar results for $d^-$ replacing 0 by 1. □

**Remark 3.1.** *From now on we use $\limsup_{n\to\infty} \frac{1}{n} S_n\phi(0)$ to refer one of the items in the above Lemma.*

The set $\mathcal{H}$ is defined as the set of $\phi \in C^\alpha([0,1], \mathcal{P})$ such that

$$P_{top}(\phi, L) > \max\left\{\limsup_{n\to\infty} \frac{1}{n}(S_n\phi)(0), \limsup_{n\to\infty} \frac{1}{n}(S_n\phi)(1)\right\}.$$

### 3.2. $\mathcal{H}$ is not empty.

To show that $\mathcal{H}$ is not empty first we recall that all one-dimensional Lorenz-like maps are topologically conjugate to a $\beta$-transformation. Let $1 < \beta \leq 2$ and $\alpha \geq 0$ such that $\alpha + \beta \leq 2$, then the map $T(x) = \beta x + \alpha \mod 1$ is called $\beta$-transformation. Viewed as a map of the interval a $\beta$-transformation has a point of discontinuity at $D = (1-\alpha)/\beta$. In addition

$$T(D^+) = 0, T(D^-) = 1, T(1) = \beta + \alpha - 1, T(0) = \alpha$$

and

$$T(x) = \begin{cases} \beta x + \alpha, & \text{if } x \in [0, (1-\alpha)/\beta) \\ \beta x + \alpha - 1, & \text{if } x \in ((1-\alpha)/\beta, 1]. \end{cases}$$

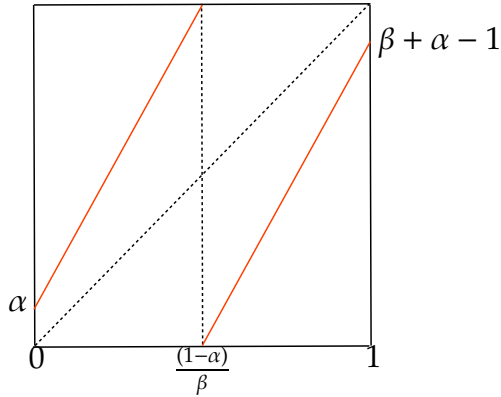

FIGURE 2. Generalized $\beta$-transformations: $T(x) = \beta x + \alpha \mod 1$.

It is known (see [6], [8], [9], [10]) that every one-dimensional expanding Lorenz-like map is topologically conjugate to a $\beta$-transformation.

In this context, Glendinning proved the following theorem:

**Theorem 3.1** ([3]). *Let $L : [0,1] \setminus \{d\} \to [0,1]$ be a one-dimensional Lorenz-like map expanding. Then $L$ is topologically conjugate to a $\beta$-transformation $T$ and $h_{top}(L) = \ln(\beta)$.*

**Proof that $\mathcal{H}$ is not empty.** Consider the potential $\overline{\phi} \equiv 0$. By Theorem 3.1, we get
$$P_{top}(\overline{\phi}, L) = h_{top}(L) = \ln(\beta) > 0,$$
since $1 < \beta \leq 2$.

On the other hand,
$$P_{top}(\overline{\phi}, L, \partial\mathcal{P}) = \max\left\{\limsup_{n\to\infty} \frac{1}{n}(S_n\overline{\phi})(0), \limsup_{n\to\infty} \frac{1}{n}(S_n\overline{\phi})(1)\right\} = 0.$$

Thus,
$$P_{top}(\overline{\phi}, L, \partial\mathcal{P}) = 0 < \ln(\beta) = h_{top}(L) = P_{top}(\overline{\phi}, L)$$
and we conclude that $\overline{\phi} \in \mathcal{H}$.

3.3. **Every element of $\mathcal{H}$ admits a unique equilibrium measure.** To show that every element of $\mathcal{H}$ admits an unique equilibrium measure we first prove that the pressure of the boundary $\partial\mathcal{P}$ can be written in terms of asymptotic values of the Birkhoff average of the boundary point of partition. We begin with the following simple result based on distortion argument.

**Lemma 3.2.** *Let $L : [0,1] \setminus \{d\} \to [0,1]$ be a one-dimensional Lorenz-like expanding map and $\phi \in C^\alpha([0,1], \mathcal{P})$. Then for $n$ large enough there exists a constant $C > 0$ such that $|(S_n\phi)(x) - (S_n\phi)(y)| \leq C$, for all $x, y \in C_n$ and $C_n \in \mathcal{P}^{(n-1)}$.*

*Proof.* If $\phi$ is $\alpha$-Hölder continuous on $\mathcal{P}$ there exists a constant $K > 0$ such that
$$|(S_n\phi)(x) - (S_n\phi)(y)| \leq \sum_{i=0}^{n-1} |\phi(L^i(x)) - \phi(L^i(y))| \leq K\sum_{i=0}^{n-1} |L^i(x) - L^i(y)|^\alpha.$$

Since $x, y \in C_n$ by property L.3 of Definition 2.1 there exists $0 < \lambda < 1$ such that
$$|(S_n\phi)(x) - (S_n\phi)(y)| \leq K\sum_{i=0}^{\infty} \lambda^{i\alpha}.$$

Hence it is enough to take $C = K\sum_{i=0}^{\infty} \lambda^{i\alpha}$.  □

**Corollary 3.1.** *The Lemma 3.2 is true on the closure of the cilinders $C_n$, denoted by $\overline{C_n}$, i.e, if $\overline{x}, \overline{y} \in \overline{C_n}$ then there $C > 0$, such that $|(S_n\phi)(\overline{x}) - (S_n\phi)(\overline{y})| \leq C$.*

*Proof.* Just write points in $\overline{C_n}$ as limits of points in $C_n$.  □

Using the Corollary 3.1, we can characterize the $P_{top}(\phi, \partial\mathcal{P}, L)$ as follows.

**Proposition 3.1.** *Let* $L : [0,1] \setminus \{d\} \to [0,1]$ *be a one-dimensional Lorenz-like expanding map and* $\phi \in C^\alpha([0,1], \mathcal{P})$, *then*

$$P_{top}(\phi, \partial \mathcal{P}, L) = \max \left\{ \limsup_{n \to \infty} \frac{1}{n}(S_n \phi)(0), \limsup_{n \to \infty} \frac{1}{n}(S_n \phi)(1) \right\}.$$

*Proof.* Consider $C_n^i$, $1 \le i \le 4$ the cylinder in $\mathcal{P}^{(n-1)}$ such that $\{0, d, 1\} \cap \partial C_n^i \ne \emptyset$, where $\partial C_n^i$ is bourdary of $C_n^i$. Since $\sup_{x \in C_n} e^{(S_n \phi)(x)} \le \sup_{x \in \overline{C_n}} e^{(S_n \phi)(x)}$ it follows that

$$P_{top}(\phi, \partial \mathcal{P}, L) \le \limsup_{n \to \infty} \frac{1}{n} \log \left( \sum_{C_n \in \mathcal{P}^{n-1} : \partial \mathcal{P} \cap \overline{C_n} \ne \emptyset} \sup_{x \in \overline{C_n}} e^{(S_n \phi)(x)} \right).$$

Using the Lemma 3.2, for $x \in C_n^i$, we have $(S_n \phi)(x) \le (S_n \phi)(b_i) + C$, where $b_i \in \{0, d^-, d^+, 1\}$ depending on whether $C_n$ is on the left or right hand side of the discontinuity $d$. It turns out that

$$P_{top}(\phi, \partial \mathcal{P}, L) \le \limsup_{n \to \infty} \frac{1}{n} \log \left( \sum_{i=1}^{4} e^{(S_n \phi)(b_i)} \right)$$

$$= \max_{1 \le j \le 4} \left\{ \limsup_{n \to \infty} \frac{1}{n} S_n \phi(b_j) \right\}$$

$$\le \max \left\{ \limsup_{n \to \infty} \frac{1}{n} S_n \phi(0), \limsup_{n \to \infty} \frac{1}{n} S_n \phi(1) \right\}$$

were in the above inequalities we used Lemma 3.1.

Conversely, since $\sup_{x \in C_n} e^{(S_n \phi)(x)} \ge \inf_{x \in \overline{C_n}} e^{(S_n \phi)(x)}$ we obtain that

$$P_{top}(\phi, \partial \mathcal{P}, L) \ge \limsup_{n \to \infty} \frac{1}{n} \log \left( \sum_{C_n \in \mathcal{P}^{n-1} : \partial \mathcal{P} \cap \overline{C_n} \ne \emptyset} \inf_{x \in \overline{C_n}} e^{(S_n \phi)(x)} \right).$$

Using again the Lemma 3.2, for $x \in C_n^i$, we have $(S_n \phi)(b_i) - C \le (S_n \phi)(x)$. Thus

$$P_{top}(\phi, \partial \mathcal{P}, L) \ge \limsup_{n \to \infty} \frac{1}{n} \log \left( e^{\max_{1 \le j \le 4} \{S_n \phi(b_j)\}} \right) \ge \max_{1 \le j \le 4} \left\{ \limsup_{n \to \infty} \frac{1}{n} S_n \phi(b_j) \right\}.$$

Applying again Lemma 3.1 we have

$$P_{top}(\phi, \partial \mathcal{P}, L) \ge \max \left\{ \limsup_{n \to \infty} \frac{1}{n} S_n \phi(0), \limsup_{n \to \infty} \frac{1}{n} S_n \phi(1) \right\}.$$

□

**Corollary 3.2.** *Let $L : [0,1]\setminus\{d\} \to [0,1]$ be a one-dimensional Lorenz-like expanding. Then the map $P_{top}(\,\cdot\,, \partial \mathcal{P}, L) : C^\alpha([0,1], \mathcal{P}) \longrightarrow \mathbb{R}$ is continuous.*

*Proof.* Note that for all $\psi, \phi \in C^\alpha([0,1], \mathcal{P})$ by Proposition 3.1 we have that

$$\left|P_{top}(\psi, \partial \mathcal{P}, L) - P_{top}(\phi, \partial \mathcal{P}, L)\right| =$$

$$= \left|\limsup_{n\to\infty} \frac{1}{n} \max\{S_n\psi(0), S_n\psi(1)\} - \limsup_{n\to\infty} \frac{1}{n} \max\{S_n\phi(0), S_n\phi(1)\}\right|$$

$$\leq 2\limsup_{n\to\infty} \frac{1}{n}\left|\frac{S_n\psi(0) - S_n\phi(0)}{2}\right| + 2\limsup_{n\to\infty} \frac{1}{n}\left|\frac{S_n\psi(1) - S_n\phi(1)}{2}\right|$$

$$\leq \|\psi - \phi\|.$$

$\square$

**Proof of that every element of $\mathcal{H}$ admits a unique equilibrium measure.** Let $\phi \in \mathcal{H}$. By definition we have that

$$\max\left\{\limsup_{n\to\infty} \frac{1}{n} S_n\phi(0), \limsup_{n\to\infty} \frac{1}{n} S_n\phi(1)\right\} < P_{top}(\phi, L).$$

Thus by proposition 3.1, we obtain

$$P_{top}(\phi, \partial \mathcal{P}, L) = \max\left\{\limsup_{n\to\infty} \frac{1}{n} S_n\phi(0), \limsup_{n\to\infty} \frac{1}{n} S_n\phi(1)\right\} < P_{top}(\phi, L).$$

On the other hand, by remark 2.1 we have that one-dimensional Lorenz-like expanding map is LEO. Thus applying theorem 2.1 we conclude that if $\phi \in \mathcal{H}$ then $\phi$ admits a unique equilibrium measure

3.4. **$\mathcal{H}$ is an open set in $C^\alpha([0,1], \mathcal{P})$.** Let us first observe that that $\mathcal{H} = \mathcal{H}_+ \cap \mathcal{H}_-$, where

$$\mathcal{H}_+ = \left\{\phi \,:\, P_{top}(\phi, L) > \limsup_{n\to\infty} \frac{1}{n}(S_n\phi)(0)\right\}$$

$$\mathcal{H}_- = \left\{\phi \,:\, P_{top}(\phi, L) > \limsup_{n\to\infty} \frac{1}{n}(S_n\phi)(1)\right\}.$$

Consider $\phi \in C^\alpha([0,1], \mathcal{P})$. To prove that $\mathcal{H}$ is $C^0$-open in $C^\alpha([0,1], \mathcal{P})$ we need the following Lemmas:

**Lemma 3.3.** $P_{top}(\cdot, L) : C^\alpha([0,1], \mathcal{P}) \to \mathbb{R}$ *defined by $\phi \longmapsto P_{top}(\phi, L)$. Then for each $\phi, \psi \in C^\alpha([0,1], \mathcal{P})$, we have $|P_{top}(\phi, L) - P_{top}(\psi, L)| \leq \|\phi - \psi\|_0$, in other words, $P_{top}(\cdot, L)$ is a Lipschitz function.*

*Proof.* The proof of the Lemma follows from Theorem 9.7 of [13], with appropriate modifications. $\square$

**Lemma 3.4.** *The functionals $P_+$, $P_-$ : $C^\alpha([0,1], \mathcal{P}) \to \mathbb{R}$, defined by*
$$\phi \longmapsto P_-(\phi) = \limsup_{n\to\infty} \frac{1}{n}(S_n\phi)(1) \text{ and } \phi \longmapsto P_+(\phi) = \limsup_{n\to\infty} \frac{1}{n}(S_n\phi)(0)$$
*are continuous.*

*Proof.* We give the proof only for the map $P_+$, the other case is analogous. As in Lemma 3.1 we consider two cases. If there is $n_0 \in \mathbb{N}$ such that $L^{n_0}(0) = d$, then the Lemma 3.4 is trivial, since $P_+$ is constant. If there is not $n_0$ such that $L^{n_0}(0) = d$, then by the definition of Birkhoff sum, we have that

$$\begin{aligned}
|P_+(\phi) - P_+(\psi)| &= \left|\limsup_{n\to\infty} \frac{1}{n}(S_n\phi)(0) - \limsup_{n\to\infty} \frac{1}{n}(S_n\psi)(0)\right| \\
&\leq \limsup_{n\to\infty} \frac{1}{n}\left|(S_n\phi(0) - S_n\psi(0))\right| \\
&\leq \|\phi - \psi\|_0.
\end{aligned}$$
(3.1)

Thus the inequality (3.1) proves the Lemma. □

**Proof of that $\mathcal{H}$ is $C^0$-open in $C^\alpha([0,1], \mathcal{P})$.** Now, let us prove the first part of the Theorem A. Note that, to prove that $\mathcal{H}$ is open in $C^\alpha([0,1], \mathcal{P})$ it is sufficient to show that $\mathcal{H}_+$ and $\mathcal{H}_-$ are open in $C^\alpha([0,1], \mathcal{P})$.

Consider $\tilde{P}$ : $\mathcal{H}_+ \to \mathbb{R}$ defined by $\tilde{P}(\phi) = P_{top}(\phi, L) - P_+(\phi)$, for all $\phi \in \mathcal{H}_+$. Combining Lemmas 3.3 and 3.4 we obtain that $\tilde{P}$ is continuous, and $\mathcal{H}_+ = \{\phi : P_{top}(\phi, L) > P_+(\phi)\} = \tilde{P}^{-1}(A)$, with $A = (0, \infty)$. Thus $\mathcal{H}_+$ is an open set of $C([0,1], \mathbb{R})$. Similarly we get $\mathcal{H}_+$ is open and we conclude that $\mathcal{H}$ is open in $C^\alpha([0,1], \mathcal{P})$.

**3.5. $\mathcal{H}$ is an dense set in $C^\alpha([0,1], \mathcal{P})$.** Fixed $\phi \in C^\alpha([0,1], \mathcal{P})$, our proof consists in the construction of a sequence of potentials $\phi_{\epsilon,k,l} \in \mathcal{H}$ such that $\phi_{\epsilon,k,l} \to \phi$ in $C^0$-topology. To this end, first by properties of $\beta$–transformation we construct an auxiliary family $\{A_n\}_{n\in\mathbb{N}}$ of subsets of $[0,1]$, generated by the partition $\mathcal{F} = \{(0, (1-\alpha)/\beta), ((1-\alpha)/\beta, 1)\}$. Then we find periodic points with large enough period inside the cylinders that are the right and left hand side of the discontinuity. Finally we construct $\phi_{\epsilon,k,l}$ by perturbing $\phi$ along the orbit of the above mentioned periodic points in order to obtain higher pressure and dominate the pressure of the boundary.

**Construction of auxiliary sets $\{A_n\}_{n\in\mathbb{N}}$.** Note that by Theorem 3.1 to find a periodic point for one-dimensional Lorenz-like map expanding is equivalent to find a periodic point for the $\beta$-transformation. Therefore, the construction of periodic points will be done using $\beta$-transformation. The boundary of such a system $([0,1], \mathcal{F}, T)$, where $\mathcal{F} = \{(0, (1-\alpha)/\beta), ((1-\alpha)/\beta, 1)\}$, is $\partial\mathcal{F} = \{0, (1-\alpha)/\beta, 1\}$.

Let $\mathcal{F}^{(n)}$ be the collection of $n$-cylinders dynamically defined by the transformation $T$, i.e., the non empty intersections $C_n = \cap_{j=1}^{n} T^{-j}(F_j)$, where $F_j \in \mathcal{F}$.

We define $C_n^+$ and $C_n^-$ being the cylinders respectively at the right and left hand side of the discontinuity $D$, i.e., $D \in \partial C_n^\pm$, where $+$ or $-$ represent the cylinders on the right or left side of $D$, respectively. We introduce an auxiliary family $A_n$ by induction as follows: Let $A_0 := ((1-\alpha)/\beta, 1)$ and for $n \geq 0$ we write

$$A_{n+1} = \begin{cases} T(A_n), & \text{if } D \notin A_n \\ T(A_n^*), & \text{if } D \in A_n, \end{cases}$$

where $A_n^*$ is the connected component of $A_n \setminus \{D\}$ containing $T^n(D^+)$.

**Definition 3.1.** *An integer $N$ is a* cutting time *for $T$ if $D \in A_N$.*

**Construction of periodic points.** In this section we find periodic points with large enough periods inside the cylinders that are on right and left hand side of the discontinuity. More precisely, we will prove using properties of $\beta$-transformations that there exists sequences of integer $\{N_k^\pm\}_{k \in \mathbb{N}}$ and a sequences periodic points such that $p_k^\pm \to d^\pm, d^\pm \in \partial C_{N_k^\pm}$ and $L^{N_k^\pm}(p_k^\pm) = p_k^\pm$ for some $p_k^\pm \in C_{N_k^\pm}$, where $+$ or $-$ represent the cylinders on the right or left side of $d$, respectively. Lastly, we show that $P(\phi, \partial \mathcal{P}, L)$ can be calculated by the average of $p_k^\pm$.

**Lemma 3.5.** *Let $N^+$ be a cutting time for $T$ and $C_{N^+} = (D^+, B^+)$ a cylinder. Suppose that $T^{N^+}(B^+) > B^+$ then there exists a periodic point $p_k^+$ of period $N^+$ for $T$ such that $p_k^+ \in C_{N^+}$. Analogous results we obtain to the cylinders at the left hand side of $D$.*

*Proof.* We first show that if $T^{N^+}(B^+) > B^+$ then there exists a periodic point of period $N^+$ for $T$ inside $C_{N^+}$. The case for $T^{N^-}(B^-) < B^-$ is similar. Observe that $T^{N^+}(C_{N^+}) = A_{N^+}$. As $N^+$ is a cutting time we obtain that $A_{N^+} = (T^{N^+}(D^+), T^{N^+}(B))$ and $D \in A_{N^+}$. So $T^{N^+}(C_{N^+}) = (T^N(D^+), T^N(B))$. As $T^{N^+}(B^+) > B^+$ and $D$ is a cutting time we obtain

$$T^{N^+}(C_{N^+}) = A_{N^+} = (T^{N^+}(D^+), T^{N^+}(B^+)) \supset (D, B^+) = C_{N^+},$$

so there exists $p_{N^+}^+ \in C_{N^+}$ such that $T^{N^+}(p_{N^+}^+) = p_{N^+}^+$. □

**Lemma 3.6.** *There exists infinity cutting times $N_0^+ < N_1^+ < \cdots < N_k^+ < \cdots$ such that $T^{N_k^+}(B^+) > B^+$, for all $k \in \mathbb{N}$. Similar results we get to the cylinders at the left hand side of $D$.*

*Proof.* We first prove that there exists infinity cutting times $N_0^+ < N_1^+ < \cdots < N_k^+ < \cdots$ such that $T^{N_k^+}(B^+) > B^+$, for all $k \in \mathbb{N}$. We can proceed similarly to prove the second case. To this end, fix $N_0^+ > 0$. By contradiction suppose that for all cutting time $N^+ \geq N_0^+$, we have that $T^{N^+}(B^+) < B^+$. Thus there exists cutting times $N_0^+ < N_1^+ < \cdots < N_k^+ < \cdots$ such that $T^{N_k^+}(B^+) < B^+$, for all $k \in \mathbb{N}$. Therefore,

$$\frac{|C_{N_{k+1}^+}|}{|C_{N_k^+}|} = \frac{|A_{N_{k+1}^+}|}{\beta^{N_{k+1}^+ - N_k^+}|A_{N_k^+}|} = \frac{|T^{N_k^+}(D^+) - D|}{|T^{N_k^+}(D^+) - T^{N_k^+}(B^+)|}. \tag{3.2}$$

The second equality above comes from the fact that there is no cutting time between $N_k^+$ and $N_{k+1}^+$ and consequently

$$|A_{N_{k+1}^+}| = \beta^{N_{k+1}^+ - N_k}|T^{N_k^+}(D^+) - D|.$$

Now,
$$\frac{|T^{N_k^+}(D^+) - D|}{|T^{N_k^+}(D^+) - T^{N_k^+}(B^+)|} \geq \frac{|T^{N_k^+}(D^+) - T^{N_k^+}(B^+)| - |T^{N_k^+}(B^+) - D|}{|T^{N_k^+}(D^+) - T^{N_k^+}(B^+)|} \tag{3.3}$$

$$= 1 - \frac{|T^{N_k^+}(B^+) - D|}{|T^{N_k^+}(D^+) - T^{N_k^+}(B^+)|} = 1 - \frac{1}{\beta^{N_k^+}}\frac{|T^{N_k^+}(B^+) - D|}{|D - B^+|}$$

Since we are considering $T^{N_k^+}(B^+) < B^+$ we have

$$-\frac{1}{\beta^{N_k^+}}\frac{|T^{N_k^+}(B^+) - D|}{|D - T(B^+)|} > -\frac{1}{\beta^{N_k^+}}. \tag{3.4}$$

Thus, using (3.2), (3.3) and (3.4) we have that

$$\frac{|C_{N_{k+1}^+}|}{|C_{N_k^+}|} > 1 - \frac{1}{\beta^{N_k^+}}.$$

As a consequence we have that

$$\frac{|C_{N_{k+1}^+}|}{|C_{N_{k_0}^+}|} = \frac{|C_{N_{k+1}^+}|}{|C_{N_k^+}|} \cdot \frac{|C_{N_k^+}|}{|C_{N_{k-1}^+}|} \cdots \cdots \frac{|C_{N_{k_0+1}^+}|}{|C_{N_{k_0}^+}|}$$

$$> \prod_{j=1}^{\infty}\left(1 - \frac{1}{\beta^j}\right) =: \gamma > 0.$$

Therefore, $|C_{N_k^+}| > \gamma \cdot |C_{N_{k_0}^+}|$, for all $k \geq k_0$. This gives us a contradiction since by construction $|C_{N_k^+}| \to 0$ when $k \to \infty$. Thus there exists infinity cutting times $N^+$ such that $T^{N^+}(B^+) > B^+$. □

**Corollary 3.3.** *Let $\mathcal{P}^{(N_k^+)}$ the collection of $N_k^+$-cylinder for one-dimensional Lorenz-like expanding map L. Then there exists $N_k^+ \in \mathbb{N}$ sucht that $d \in \partial C_{N_k^+}$ and $L^{N_k^+}(p_k^+) = p_k^+$ for some $p_k^+ \in C_{N_k^+}$. To the cylinder $C_{N_k^-}$ the construction is analogous.*

*Proof.* Since one-dimensional Lorenz-like expanding map $L$ is topologically conjugate the maps $T$ it follows immediately from Lemma 3.5 that periodic point $p_k^\pm$ there exists. □

**Lemma 3.7.** *Let $L : [0,1] \setminus \{d\} \to [0,1]$ be a one-dimensional Lorenz-like expanding map and consider $\phi \in C^\alpha([0,1], \mathcal{P})$. If $N_k^+ \in \mathbb{N}$ such that $d \in \partial C_{N_k^+}$, $p_k^+ \in C_{N_k^+}$ and $L^{N_k^+}(p_k^+) = p_k^+$, then*

$$(3.5) \qquad \limsup_{n \to \infty} \left( \lim_{k \to \infty} \frac{1}{n} S_n \phi(p_k^+) \right) = \limsup_{k \to \infty} \left( \lim_{n \to \infty} \frac{1}{n} S_n \phi(p_k^+) \right),$$

*Replacing $p_k^+$ by $p_k^-$ we get similar results.*

*Proof.* We give the proof only for $p_k^+$. The case $p_k^-$ is similar. To see this, we first observe that as $\sup(\phi) < \infty$ and then

$$\left| \frac{1}{n} S_n \phi(p_k^+) \right| \leq \frac{1}{n} \sum_{j=0}^{n-1} |\phi(L^j(p_k^+))| \leq \sup(\phi), \text{ for all } n, k \in \mathbb{N}.$$

Thus the double superior limit $\limsup_{n,k \to \infty} \frac{1}{n} S_n \phi(p_k^+)$ exists. Now, for fixed $k \in \mathbb{N}$, we can be write

$$n = q_n^+ N_k^+ + r_n^+, \ 0 \leq r_n^+ \leq N_k^+ - 1.$$

Thus we have that the following limit exists

$$(3.6) \quad \begin{aligned} \lim_{n \to \infty} \frac{1}{n} S_n \phi(p_{k^+}) &= \lim_{n \to \infty} \left( \frac{q_n^+}{q_n^+ N_k^+ + r_n^+} S_{N_k^+} \phi(p_k^+) + \frac{1}{q_n^+ N_k^+ + r_n^+} S_{r_n^+} \phi(p_k^+) \right) \\ &= \frac{1}{N_k^+} S_{N_k^+} \phi(p_k^+), \end{aligned}$$

where we get the last equality because, as $0 \leq r_n^+ \leq N_k^+ - 1$, then $S_{r_n^+} \phi(p_k^+)$ is limited. On the other hand, for fixed $n \in \mathbb{N}$, as $p_k^+ \to d^+$ the following limit exists by continuity, *i.e*,

$$(3.7) \qquad \lim_{k \to \infty} \frac{1}{n} S_n \phi(p_k^+) = \frac{1}{n} S_n \phi(d^+).$$

Combining (3.6) and (3.7) we see at once that equal (3.5) is true.

□

**Corollary 3.4.** *Let $L : [0,1] \setminus \{d\} \to [0,1]$ be a one-dimensional Lorenz-like expanding map and consider $\phi \in C^\alpha([0,1], \mathcal{P})$. If $p_k^+$ is sequence of periodic points satisfying the conditions Lemma 3.7 then*

$$\limsup_{k \to \infty} \frac{1}{N_k^+} S_{N_k^+} \phi(p_k^+) = \limsup_{n \to \infty} \frac{1}{n} S_n \phi(0).$$

*Replacing* 0 *by* 1 *and* $p_k^+$ *by* $p_k^-$ *we get similar results.*

*Proof.* For each $n$, as $\lim_{k \to \infty} p_k^+ = d^+$, we obtain

$$\text{(3.8)} \qquad \lim_{k \to \infty} \frac{1}{n} S_n \phi(p_k) = \frac{1}{n} S_n \phi(d^+).$$

Letting $n \to \infty$ we can rewrite (3.8) as

$$\limsup_{n \to \infty} \left( \lim_{k \to \infty} \frac{1}{n} S_n \phi(p_k) \right) = \limsup_{n \to \infty} \frac{1}{n} S_n \phi(d^+).$$

By Lemma 3.7 we obtain

$$\limsup_{k \to \infty} \left( \lim_{n \to \infty} \frac{1}{n} S_n \phi(p_k) \right) = \limsup_{n \to \infty} \frac{1}{n} S_n \phi(d^+).$$

As

$$\lim_{n \to \infty} \frac{1}{n} S_n \phi(p_k^+) = \frac{1}{N_k} S_{N_k} \phi(p_k)$$

this finishes the proof. $\square$

**Construction of potential $\phi_{\epsilon,k,l}$.** Recall that by Corollary 3.4, there exists subsequence $N_k^+ \to \infty$ such that $d \in \partial C_{N_k^+}$ and $p_k^\pm \in C_{N_k^\pm}$ such that $L^{N_k^\pm}(p_k) = p_k^\pm$. Let $I_j^\pm = (L^j(p_k^\pm) - \delta_k, L^j(p_k^\pm) + \delta_k)$ be intervals, where $0 \leq j \leq N_k^\pm - 1$. Since the orbit $p_k^\pm$ is finite, there exists $\delta_k^\pm > 0$ such that $I_i^\pm \cap I_j^\pm = \emptyset$, for all $0 \leq j, i \leq N_k^\pm - 1$ with $i \neq j$. Here and subsequently, we denote $\delta_{k,l} = \dfrac{\delta_k}{l}$. Let

$$B_{\epsilon,k,l}^\pm(x) = \begin{cases} \displaystyle\sum_{j=0}^{N_k^\pm - 1} B_{\epsilon,k,j,l}^\pm(x) & , \ x \in \displaystyle\bigcup_{j=0}^{N_k - 1} I_j^\pm \\ 0 & , \ \text{otherwise.} \end{cases}$$

where $B_{\epsilon,k,j,l}^\pm$ is bump function defined by figure 3.5.

**Lemma 3.8.** *Let* $\phi_{\epsilon,k,l}^\pm(x) = \phi(x) + B_{\epsilon,k,l}^\pm(x)$, *where* $\phi \in C^\alpha([0,1], \mathcal{P})$. *Then we have the following properties:*

(1) $\phi_{\epsilon,k,l}^\pm$ *is Holder continuous;*
(2) $\|\phi_{\epsilon,k,l}^\pm - \phi\|_{C^0} < \epsilon$, *for all* $\phi \in C^\alpha([0,1], \mathcal{P})$;
(3) *As* $\phi_{\epsilon,k,l}^\pm$ *is built on the orbit of the periodic point* $p_k^+$, *we have that*

$$S_{N_k} \phi_{\epsilon,k,l}^\pm(p_k^\pm) = S_{N_k} \phi(p_k^\pm) + \epsilon.$$

*Proof.* The proof follows of the construction of potential $\phi_{\epsilon,k,l}^\pm$. $\square$

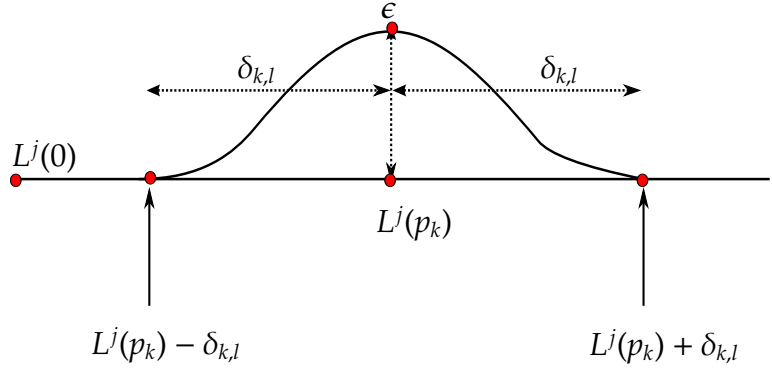

FIGURE 3. Bump function $B^{\pm}_{\epsilon,k,j,l}$.

**Lemma 3.9.** *Let $O(p^{\pm}_k)$ the orbit of $p^{\pm}_k$ and $\chi_{O(p^{\pm}_k)}$ the characteristic function of $O(p^{\pm}_k)$. Then $\lim_{l\to\infty} B^{\pm}_{\epsilon,k,l}(x) = \epsilon \cdot \chi_{O(p^{\pm}_k)}(x)$, for all $x \in [0,1]$.*

*Proof.* If $x \in O(p^{\pm}_k)$, then there exists $j_0 \in \{0, 1, \cdots, N^{\pm}_k - 1\}$ such that $x = L^{j_0}(p^{\pm}_k)$. Then $B^{\pm}_{\epsilon,k,l}(x) = B^{\pm}_{\epsilon,k,l}(L^{j_0}(p^{\pm}_k)) = \sum_{j=0}^{N^{\pm}_k-1} B^{\pm}_{\epsilon,k,l,j}(L^{j_0}(p^{\pm}_k)) = \epsilon$, since by the construction $B^{\pm}_{\epsilon,k,l,j}(L^{j_0}(p^{\pm}_k)) = 0$, whenever $j \neq j_0$, and $B^{\pm}_{\epsilon,k,l,j_0}(L^{j_0}(p^{\pm}_k)) = \epsilon$. Therefore, we conclude that $\lim_{l\to\infty} B^{\pm}_{\epsilon,k,l}(x) = \epsilon$.

On the other hand, let $x \in [0,1] \setminus O(p^{\pm}_k)$. If $x \notin \bigcup_{j=0}^{N^+_{\pm}} I^{\pm}_j$ then $\lim_{l\to\infty} B^{\pm}_{\epsilon,k,l}(x) = 0$. Otherwise if $x \in I^{\pm}_j$ for some $j \in \{0, 1, \cdots, N^{\pm}_k - 1\}$. However, since $\delta_{k,l} \to 0$ as $l \to \infty$, then for large enough $l_0$ we also have $x \notin \bigcup_{j=0}^{N^+_{\pm}} I^{\pm}_j$ and $\lim_{l\to\infty} B^{\pm}_{\epsilon,k,l}(x) = 0$. Therefore, for all $x \in [0,1]/O(p^{\pm}_k)$, we get that $\lim_{l\to\infty} B^{\pm}_{\epsilon,k,l}(x) = 0$ this complete the proof. □

**Lemma 3.10.** *The pressure the of bump function $B^{\pm}_{\epsilon,k,l}$ is null, i.e.,*

$$\lim_{l\to\infty} P_{top}(B^{\pm}_{\epsilon,k,l}, \partial \mathcal{P}, L) = 0.$$

*Proof.* By the Lemma 3.9, we have that $\lim_{l\to\infty} B^{\pm}_{\epsilon,k,l}(x) = \epsilon \cdot \chi_{O(p^{\pm}_k)}(x)$. As $O(p^{\pm}_k)$ is the closed set we obtain that $\chi_{O(p^{\pm}_k)} : [0,1] \to \mathbb{R}$ is upper semi-continuous. Thus, applying Corollary 3.2 we have that

$$\lim_{l \to \infty} P_{top}(B^{\pm}_{\epsilon,k,l}, \partial \mathcal{P}, L) = P_{top}(\epsilon \cdot \chi_{O(p_k^{\pm})}, \partial \mathcal{P}, L)$$

$$= \epsilon \cdot \left( \limsup_{n \to \infty} \frac{1}{n} \max \left\{ S_n \chi_{O(p_k^{\pm})}(0), S_n \chi_{O(p_k^{\pm})}(1) \right\} \right)$$

$$= 0,$$

as $L^j(x) \notin O(p_k^+)$ and $\chi_{O(p_k^+)}(L^j(x)) = 0$, for all $0 \leq j \leq n-1$ and this finishes the proof. $\square$

**Corollary 3.5.** *Define* $\phi_{\epsilon,k,l}(x) = max\{\phi^+_{\epsilon,k,l}(x), \phi^-_{\epsilon,k,l}(x)\}$, *for all* $x \in [0,1]$. *Then*

$$\lim_{l \to \infty} P_{top}(\phi_{\epsilon,k,l}, \partial \mathcal{P}, L) = P_{top}(\phi, \partial \mathcal{P}, L).$$

**Proof of that $\mathcal{H}$ is $C^0$-dense in $C^\alpha([0,1], \mathcal{P})$.** Our proof starts by recalling the follows result proved by Buzzi and Sarig in [2], where $L$ is the Lorenz map.

**Proposition 3.2** ([2]). *Consider $\phi$ a piecewise uniformly continuous potential and let $\nu$ be an ergodic probability measure. If $\nu(S) > 0$, then $P_{top}(\phi, S, L) \geq h_\nu(L) + \int \phi \, d\nu$, where $h_\nu(L)$ is the metric entropy of $\nu$.*

Let $\phi \in C^\alpha([0,1], \mathcal{P})$. Remember that our proof consists in the construction of a potential $\phi_{\epsilon,k,l}$ such that $\|\phi_{\epsilon,k,l} - \phi\|_{C^0} < \epsilon$ and

(3.9) $$P_{top}(\phi_{\epsilon,k,l}, L) > P_{top}(\phi_{\epsilon,k,l}, \partial \mathcal{P}, L),$$

i.e, $\phi_{\epsilon,k,l} \in \mathcal{H}$. To this end, fix $\epsilon > 0$ and consider

$$\phi_{\epsilon,k,l}(x) = max\{\phi^+_{\epsilon,k,l}(x), \phi^-_{\epsilon,k,l}(x)\}, \text{ for all } x \in [0,1].$$

By the Lemma 3.8 the potential $\phi_{\epsilon,k,l}$ satisfies the condition $\|\phi_{\epsilon,k,l} - \phi\|_{C^0} < \epsilon$, for all $k, l$. Thus we only need to show that the inequality (3.9) is true. By definition of $P_{top}(\phi_{\epsilon,k,l}, L)$, we always have that there exists $k_0, l_0$ sucht that

$$P_{top}(\phi_{\epsilon,k,l}, L) \geq P_{top}(\phi_{\epsilon,k,l}, \partial \mathcal{P}, L)$$

for all $k \geq k_0$ and $l \geq l_0$. To obtain a contradiction, suppose that for all $k, l$,

(3.10) $$P_{top}(\phi_{\epsilon,k,l}, L) = P_{top}(\phi_{\epsilon,k,l}, \partial \mathcal{P}, L).$$

Consider $C_{N_k^+} \in \mathcal{P}^{(N_k^+ - 1)}$ as in Proposition 3.5, then there exists $p_k^+ \in C_{N_k^+}$ such that $L^{N_k^+}(p_k^+) = p_k^+$. Furthermore one can construct a measure $\mu_k^+(\cdot) = \left( \frac{1}{N_k^+} \sum_{j=0}^{N_k^+ - 1} \delta_{L^j(p_k^+)} \right)(\cdot)$, where $\delta_{L^j(p_k^+)}$ is the Dirac measure with $\delta_{L^j(p_k^+)}(L^j(p_k^+)) = 1$,

$j \in \{0, 1, \cdots, N_k - 1\}$. As $\mu_k^+(C_{N_k}) > 0$ by Proposition 3.2 we have

$$P_{top}(\phi_{\epsilon,k,l}, L) \geq P_{top}(\phi_{\epsilon,k,l}^+, L) \geq P_{top}(\phi_{\epsilon,k,l}^+, C_{N_k}, L)$$

$$\geq h_{\mu_k^+}(L) + \int \phi_{\epsilon,k,l}^+ d\mu_k^+$$

$$= \frac{1}{N_k} \sum_{j=0}^{N_k-1} \phi_{\epsilon,k,l}^+(L^j(p_k^+)) = \frac{1}{N_k}(S_{N_k}\phi_{\epsilon,k,l}^+)(p_k^+)$$

$$= \frac{1}{N_k}(S_{N_k}\phi)(p_k^+) + \epsilon,$$

where in the last equality we are applying the Lemma 3.8. We get the same conclusion for $p_k^- \in C_{N_k^-}$, i.e.,

$$P_{top}(\phi_{\epsilon,k,l}, L) \geq \frac{1}{N_k^-}(S_{N_k^-}\phi)(p_k^-) + \epsilon.$$

Thus,

$$(3.11) \quad P_{top}(\phi_{\epsilon,k,l}, L) \geq \max\left\{\frac{1}{N_k^+}(S_{N_k^+}\phi)(p_k^+), \frac{1}{N_k^-}(S_{N_k^-}\phi)(p_k^-)\right\} + \epsilon.$$

Letting $l \to \infty$ in the inequation (3.11) and combining (3.10) with Corollary 3.5 we get

$$P_{top}(\phi, \partial \mathcal{P}, L) \geq \max\left\{\frac{1}{N_k^+}(S_{N_k^+}\phi)(p_k^+), \frac{1}{N_k^-}(S_{N_k^-}\phi)(p_k^-)\right\} + \epsilon.$$

Indeed letting $k \to \infty$ and combining Corollary 3.4 with Proposition 3.1 we obtain

$$P_{top}(\phi, \partial \mathcal{P}, L) \geq \max\left\{\limsup_{k \to \infty} \frac{1}{N_k^+}(S_{N_k^+}\phi)(p_k^+), \limsup_{t \to \infty} \frac{1}{N_k^-}(S_{N_k^-}\phi)(p_k^-)\right\} + \epsilon$$

$$= \max\left\{\limsup_{n \to \infty} \frac{1}{n}(S_n\phi)(0), \limsup_{n \to \infty} \frac{1}{n}(S_n\phi)(1)\right\} + \epsilon$$

$$= P_{top}(\phi, \partial \mathcal{P}, L) + \epsilon.$$

that contradicts (3.10). Thus we show that the inequality (3.9) is true, which completes the proof of the Theorem A.

Faculadade de Matemática, FAMAT-UFU Uberlândia-MG, Brazil
E-mail address: jgoler@ufu.br

Faculadade de Matemática, FAMAT-UFU Uberlândia-MG, Brazil
E-mail address: bronzi@ufu.br